\ifdefined\pdfoutput
\pdfoutput=1  %% try to force PDFLaTeX processing by arXiv
\fi

\documentclass{arxiv-en-2021}

\newcommand{\ie}{i.e.{}}

\providecommand{\End}{\operatorname{End}}
\providecommand{\id}{\mathrm{id}}

\title[Proof of Cayley-Hamilton theorem]{%
  Proof of Cayley-Hamilton theorem\\
  using polynomials over the algebra\\
  of module endomorphisms}

\author{Alexey Muranov}

\date{\today}

\address{I2M, UMR 7373, CNRS, Université d'Aix-Marseille\\
  Marseille, France}

\email{\href{mailto:alexey.muranov@univ-amu.fr}
  {\nolinkurl{alexey.muranov@univ-amu.fr}}}



\begin{document}

\begin{abstract}
  If $R$ is a commutative unital ring and $M$ is a unital $R$-module,
  then each element of $\End_R(M)$ determines a left
  $\End_{R}(M)[X]$-module structure on $\End_{R}(M)$,
  where $\End_{R}(M)$ is the $R$-algebra of endomorphisms of $M$
  and $\End_{R}(M)[X] =\End_{R}(M)\otimes_RR[X]$.
  These structures provide a very short proof of the Cayley-Hamilton
  theorem, which may be viewed as a reformulation of the proof in
  \emph{Algebra} by Serge Lang.
  Some generalisations of the Cayley-Hamilton theorem can be easily
  proved using the proposed method.
\end{abstract}

\maketitle

\section{Introduction}

\begin{theorem}[Cayley-Hamilton theorem]
  Let\/ $R$ be a commutative unital ring\/
  \textup(\ie\textup, with $1_R$\/\textup)
  and\/ $M$ be a finite-rank free unital\/ $R$-module
  \textup(\ie\textup, which respects $1_R$\/\textup)\textup.
  Let\/ $a\colon M\to M$ be an endomorphism of\/ $M$ and\/
  $\chi_a\in R[X]$ be the characteristic polynomial of $a$\textup.
  Then\/ $\chi_a(a) = 0$ \textup(in\/ $R[a]\subset\End_R(M)$\textup)\textup.
\end{theorem}

The goal of this note is to show how basic properties of tensor products
provide a very short proof of this theorem and allow to generalise it.

The two main ingredients of the proof presented in this note are:
\begin{enumerate}
  \item
    the canonical isomorphism $\End_{R}(M)[X]\cong\End_{R[X]}(M[X])$,
  \item
    certain left actions of $\End_R(M)[X]$ on $\End_R(M)$
    associated to elements of $\End_{R}(M)$.
\end{enumerate}

The presented proof is essentially a reformulation of the one in
\emph{Algebra} by Serge Lang\footnotemark{} \cite{Lang.2002.a},
eliminating the need to work with bases or with matrices explicitly.
However, the author is unaware of the considered actions' of $\End_R(M)[X]$
on $\End_R(M)$ having been used in the literature before to prove
the Cayley-Hamilton theorem or generalisations thereof.

\footnotetext{%
  The proof found in the current version of
  \href{https://ncatlab.org/nlab/show/characteristic+polynomial}%
  {\emph{characteristic polynomial}}
  web page of \emph{nLab} wiki \cite{nLab.2021}
  follows the one from the Lang's book.
}

The proof by Bourbaki in \emph{Algèbre} \cite{Bourbaki.2007.edm-a-c13.fr}
is essentially different and more involved.
Not only they work with matrices explicitly, but they
also need to prove the identity $\tilde aa = a\tilde a$.

\section{Basic definitions and properties}

Necessary definitions and basic properties of modules,
their tensor products, and their exterior powers
may be found, for example, in expository papers by
Keith Conrad~\cite{Keith_Conrad.2021.expository_papers}.

Let $R$ be a commutative unital ring and $M$ be a free unital
$R$-module of finite rank~$n$.

The following usual notation shall be used:
$R[X]$ is the ring of polynomials in $X$ over $R$,
$M[X] = M\otimes_RR[X]$,
$\End_{R}(M)[X] =\End_{R}(M)\otimes_RR[X]$.

Following a common practice, elements of $R$ may be viewed as elements of
$R[X]$ or as elements of $\End_{R}(M)$ (as scalar endomorphisms),
elements of $M$ may be viewed as elements of $M[X]$, etc.

Since $M$ is free of finite rank, there is a canonical isomorphism
\[
  \End_{R}(M)[X]\cong\End_{R[X]}(M[X]).
\]
Using this isomorphism, elements of $\End_{R}(M)[X]$ may be viewed as
elements of $\End_{R[X]}(M[X])$ and vice versa.

For an endomorphism $a$ of $M$, the \emph{determinant} of $a$ is defined
by the identity
\[
  ax_1\wedge\dotsb\wedge ax_n
  = (\det a)(x_1\wedge\dotsb\wedge x_n)\qquad
  (x_1,\dotsc,x_n\in M).
\]
The \emph{adjugate endomorphism} $\tilde a$ of $a$ is defined
by the identity
\[
  ax_1\wedge\dotsb\wedge ax_{n-1}\wedge y
  = x_1\wedge\dotsb\wedge x_{n-1}\wedge\tilde ay\qquad
  (x_1,\dotsc,x_{n-1}, y\in M).
\]
Replacing $y$ with $ax_n$ in the last identity, it can be deduced that
\[
  \tilde aa = (\det a)\id_M =\det a
\]
(identifying scalar endomorphisms of $M$ with elements of~$R$).

The \emph{characteristic polynomial} of $a\in\End_R(M)$ is
the polynomial $\chi_a\in R[X]$ defined as\footnotemark{}
\[
  \chi_a =\det(a - X)
\]
(where $a - X\in\End_{R[X]}(M[X])\cong\End_{R}(M)[X]$).

\footnotetext{%
  The term \emph{characteristic polynomial} is alternatively
  (and possibly more commonly) used to denote $\det(X - a)$.
}

It is not hard to show that the degree of $\chi_a$ is $n$,
and that its leading coefficient is~$(-1)^n$.
These facts shall not be used in this note however.

Denote
\[
  t_a = a - X.
\]
Then
\[
  \chi_a =\det t_a =\tilde{t_a}t_a.
\]

\section{Proof of Cayley-Hamilton theorem through an action%
  \texorpdfstring{\\}{}
  of \texorpdfstring{$\End_R(M)[X]$}{End\_R(M)[X]}
  on \texorpdfstring{$\End_R(M)$}{End\_R(M)}}

Given $a\in\End_R(M)$,
consider the left action of $\End_R(M)[X]$ on
the $R$-module $\End_R(M)$ (forgetting its algebra structure)
denoted by the binary operator ``$\triangleleft_a$''
and defined by the rules:
\[
  f\triangleleft_ag = fg\quad\text{for}\quad f\in\End_R(M),
  \quad\text{and}\quad
  X\triangleleft_ag = ga,
\]
where $g$ is an arbitrary element of $\End_R(M)$ acted upon.
Thus, if
\[
  p = f_0 + f_1X +\dotsb + f_kX^k\in\End_R(M)[X]
  \quad\text{and}\quad
  g\in\End_R(M),
\]
then
\[
  p\triangleleft_ag = f_0g + f_1ga +\dotsb + f_kga^k\in\End_R(M),
\]
% \newpage  %% <<<<<<<<<<<<<<<<<<<<<<<<<<<<<<<<<<<<<<<<<<<<<<<<<<<<<<<<<<<<<<<
and, in particular,\footnotemark{}
\[
  p\triangleleft_a\id_M = f_0 + f_1a +\dotsb + f_ka^k\in\End_R(M).
\]

\footnotetext{%
  Since $\End_R(M)$ in general is not commutative, the element
  $f_0 + f_1a +\dotsb + f_ka^k$ of $\End_R(M)$ should not be viewed as
  the result of ``substitution'' of $a$ for $X$ in $p$.
  It may be viewed though as the result of \emph{right substitution},
  and $f_0 + af_1 +\dotsb + a^kf_k$ may be viewed as the result of
  \emph{left substitution}.
}

Thus,
\[
  (a - X)\triangleleft_a\id_M = a - a = 0
  \quad\text{and}\quad
  \chi_a\triangleleft_a\id_M = \chi_a(a).
\]

\begin{proof}[Proof of Cayley-Hamilton theorem]
  \[
    \chi_a(a) =\chi_a\triangleleft_a\id_M
    = (\tilde{t_a}t_a)\triangleleft_a\id_M
    =\tilde{t_a}\triangleleft_a(t_a\triangleleft_a\id_M)
    =\tilde{t_a}\triangleleft_a0
    = 0.
    \qedhere
  \]
\end{proof}

\section{Generalisation}

The method used above to prove the Cayley-Hamilton theorem allows to prove
seemingly more general statements, such as the following one.

\begin{proposition}
  Let\/ $R$ and\/ $M$ be as before\textup.
  Let\/ $f_1,\dotsc,f_n,a_1,\dotsc,a_n$ be endomorphisms
  of $M$ such that\/\textup:
  \begin{enumerate}
    \item
      $f_1a_1 +\dotsb + f_na_n = 0$\textup,
    \item
      $a_1,\dotsc,a_n$ commute pairwise\textup.
  \end{enumerate}
  Let
  \[
    p = f_1X_1 +\dotsb + f_nX_n
    \in\End_{R[X_1,\dotsc,X_n]}(M[X_1,\dotsc,X_n])
  \]
  and
  \[
    P =\det p\in R[X_1,\dotsc,X_n].
  \]
  Then
  \[
    P(a_1,\dotsc,a_n) = 0.
  \]
\end{proposition}

\begin{example}
  Let\/ $a$ and\/ $b$ be two commuting endomorphisms of\/ $M$\textup,
  and let
  \[
    p = bX - aY\in\End_{R[X,Y]}(M[X,Y]).
  \]
  Then substituting\/ $a$ for\/ $X$ and\/ $b$ for\/ $Y$ in\/
  $P =\det p$ yields\/~$0$\textup:
  \[
    P(a, b) = 0.
  \]
\end{example}

\printbibliography

@Book{Bourbaki.2007.edm-a-c13.fr,
  author    = {Nicolas {Bourbaki}},
  title     = {{Alg\`ebre}},
  edition   = {Reprint of the 1970 original},
  isbn      = {3-540-33849-7/pbk},
  language  = {French},
  pages     = {xiii + 636},
  publisher = {Berlin: Springer},
  subtitle  = {{Chapitres 1 \`a 3}},
  url       = {https://www.springer.com/gp/book/9783540338499},
  msc2010   = {00A05},
  year      = {2007},
  zbl       = {1111.00001},
}

@Book{Lang.2002.a,
  author    = {Serge {Lang}},
  date      = {2002},
  title     = {{Algebra}},
  edition   = {Revised third edition},
  isbn      = {0-387-95385-X/hbk},
  pages     = {xv + 918},
  publisher = {New York, NY: Springer},
  series    = {Graduate Texts in Mathematics},
  subtitle  = {{Volume 1}},
  url       = {https://www.springer.com/gp/book/9780387953854},
  volume    = {211},
  issn      = {0072-5285; 2197-5612/e},
  msc2010   = {00A05 12-01 13-01 15-01 16-01 18-01 20-01 14-01 11-01},
  zbl       = {0984.00001},
}

@Online{Keith_Conrad.2021.expository_papers,
  author  = {Keith {Conrad}},
  title   = {Expository papers},
  url     = {https://kconrad.math.uconn.edu/blurbs/},
  urldate = {2021-04-06},
}

@Online{nLab.2021,
  title      = {{nLab}},
  titleaddon = {A wiki devoted to Mathematics, Physics, and Philosophy},
  url        = {https://ncatlab.org/},
  urldate    = {2021-04-06},
}

\end{document}